\documentclass[12pt]{article}
\usepackage{amssymb,latexsym}
\def\RR{{\mathbb R}}

\newtheorem{formula}{}[section]

\newtheorem{definition}[formula]{\indent Definition}
\newtheorem{corollary}[formula]{\indent Corollary}
\newtheorem{remark}[formula]{\indent Remark}
\newtheorem{lemma}[formula]{\indent Lemma}
\newtheorem{theorem}[formula]{\indent Theorem}

\def\thrm{\begin{theorem}}
\def\thrml#1{\begin{theorem}\label{#1}}
\def\ethrm{\end{theorem}}
\def\rmrk{\begin{remark}}
\def\rmrkl#1{\begin{remark}\label{#1}}
\def\ermrk{\end{remark}}
\def\dfntn{\begin{definition}}
\def\dfntnl#1{\begin{definition}\label{#1}}
\def\edfntn{\end{definition}}
\def\nmrt{\begin{enumerate}}
\def\enmrt{\end{enumerate}}

\def\qtn{\begin{equation}}
\def\qtnl#1{\begin{equation}\label{#1}}
\def\eqtn{\end{equation}}
\def\lmm{\begin{lemma}}
\def\lmml#1{\begin{lemma}\label{#1}}
\def\elmm{\end{lemma}}
\def\crllr{\begin{corollary}}
\def\crllrl#1{\begin{corollary}\label{#1}}
\def\ecrllr{\end{corollary}}

\begin{document}
\title{}
\date{}
\maketitle
%\nopagenumbers
%\begin{titlepage}
%\title
\vspace{-0,1cm} \centerline{\bf Entropy of tropical holonomic sequences}
%\centerline{\bf AND RELATIVE KOLMOGOROV COMPLEXITY}
\vspace{7mm}
\author{
\centerline{Dima Grigoriev}
%\\[-1pt]
\vspace{3mm}
%\centerline{$^1$ St.Petersbourg University,
% Universitetskaya nab., 7/9,}
%\centerline{ St.Petersbourg,
% 199164, RUSSIA}
%\vspace{3mm}
%Dima Grigoriev \\[-1pt]
\centerline{CNRS, Math\'ematique, Universit\'e de Lille, Villeneuve
d'Ascq, 59655, France} \vspace{1mm} \centerline{e-mail:\
dmitry.grigoryev@univ-lille.fr } \vspace{1mm}
\centerline{URL:\ http://en.wikipedia.org/wiki/Dima\_Grigoriev} }
%\date{}
%\maketitle

\begin{abstract}
We introduce tropical holonomic sequences of a given order and calculate their entropy in case of the second order.
\end{abstract}

{\bf keywords}: tropical holonomic sequences, entropy


\section*{Introduction}
We consider a tropical analog of univariate holonomic sequences \cite{K}, \cite{Z}. Necessary concepts of tropical algebra
one can find in \cite{MS}.

Let $A_0(i),A_1(i),\dots,A_n(i)\in \RR[i]$ be polynomials. We say that a {\it sequence $w:=\{w_i\in \RR\, :\, i\ge 0\}$
satisfies vector $(A_0,\dots,A_n)$} if the following tropical equation
\begin{eqnarray}\label{-1}
\min \{w_i+A_0(i),\, w_{i+1}+A_1(i),\dots,w_{i+n}+A_n(i)\}
\end{eqnarray} 
holds for any $i\ge 0$. According to tropical algebra \cite{MS} this means that the minimum in (\ref{-1}) is attained for at least
two different $j$ among $w_{i+j}+A_j(i),\, 0\le j\le n$. In this case we say that $w$ is a a {\it tropical holonomic sequence of order $n$}.

For a classical holonomic sequence which satisfies equations $\sum_{0\le j\le n} A_j(i)w_{j+i}=0$ its element $w_s$ is determined uniquely for $s$ greater than $n-1$ and greater than the roots of $A_n$.

This is not the case in the tropical setting. Therefore, in \cite{G} we consider a set $W_N\subset \RR^N$ consisting of sequences $\{w_i\, : \, 0\le i<N\}$ satisfying (\ref{-1}). Then $W_N$ is a tropical prevariety \cite{MS}, so a union of a finite number of convex polyhedra. 
We define the {\it tropical entropy of vector $A_0,\dots,A_n$} as
\begin{eqnarray}\label{-2}
H:=H(A_0,\dots,A_n):=\lim \sup_{N\to \infty} \dim(W_N)/N
\end{eqnarray}
Evidently, $0\le H\le 1$. Moreover, the proof of Proposition 5.1 \cite{G} can be literally transfered to conclude that $H\le 1-1/n$.
%Moreover, if $N=p+q$ and $f:\RR^N\to \RR^p$ (respectively, $g:\RR^N \to \RR^q$) is the projection onto the first $p$ (respectively, the last $q$) coordinates then $f(W_N)\subset W_p,\, g(W_N)\subset W_q$, hence $\dim(W_N)\le \dim(W_p)+\dim(W_q)$. Therefore, due to  Fekete's subadditive lemma \cite{S} there exists a limit (which coincides as well with the infimum of this sequence)
%\begin{eqnarray}\label{-2}
%H:=H(A_0,\dots,A_n):=\lim_{N\to \infty} \dim(W_N)/N
%\end{eqnarray}
%which we call the {\it tropical entropy of vector $A_0,\dots,A_n$}. Evidently, $0\le H\le 1$. Moreover, the proof of Proposition 5.1 \cite{G} can be literally transfered to conclude that $H\le 1-1/n$.

In \cite{G} we studied the case of constant polynomials $A_0,\dots,A_n$. In this case the limit $\lim_{N\to \infty} \dim(W_N)/N (=H)$ does exist due to  Fekete's subadditive lemma \cite{S}.
In addition, one can consider Newton's graph $P$ consisting of the points $(i,\, A_i),\, 0\le i\le n$ on the plane. It was proved in Corollary 5.7 \cite{G}  that the entropy $H$ vanishes  iff $P$ is convex and $(i,\, A_i),\, 0\le i\le n$ form its vertices. Moreover, Corollary 5.7 \cite{G} states that if $H>0$ then $H\ge 1/6$. 

In the present paper in case of the second order $n=2$ we prove the existence of the limit $\lim_{N\to \infty} \dim(W_N)/N (=H)$    and explicitly calculate the tropical entropy $H$ (\ref{-2}). It appears that in this case the entropy is either $0$, either $1/4$ or $1/3$.

\section{Entropy of tropical holonomic sequences of the second order} 

Let $A(i),\, B(i),\, C(i)\in \RR[i]$ and $\{u_i\, :\, i\ge 0\}$ be a tropical holonomic sequence (of the second order) satisfying vector $(A,\, B,\, C)$, i.~e. 
\begin{eqnarray}\label{0}
\min\{u_j+A(j),\, u_{j+1}+B(j),\, u_{j+2}+C(j)\}
\end{eqnarray}
for $j\ge 0$. Denote the entropy $H:=H(A,\, B,\, C)$ (see (\ref{-2})). We write for a polynomial $A>0$ if $A(j)>0$ for $j>>0$.

\begin{theorem}
The limit $\lim_{N\to \infty} \dim(W_N)/N (=H)$ does exist and 

1) if
\begin{eqnarray}\label{1}
A(i)+C(i-1)=B(i-1)+B(i)
\end{eqnarray}
for all $i$ then $H=1/3$;

2) if 
\begin{eqnarray}\label{2}
A(i)+C(i-1)<B(i-1)+B(i),
\end{eqnarray}
\begin{eqnarray}\label{3}
A(i+1)-A(i)+C(i)-C(i-1)=B(i+1)-B(i-1)
\end{eqnarray}
then $H=1/4$:

3) otherwise $H=0$.
\end{theorem} 

In the next two sections we prove the Theorem.

\section{Case of a positive entropy}

\subsection{Case $A(i)+C(i-1)=B(i-1)+B(i)$}\label{1.1}

Under the assumption (\ref{1}) we construct a family of holonomic sequences $\{u_i\, :\, i\ge 0\}$ satisfying vector $(A,\, B,\, C)$ by induction on $i$.
As an induction base put $u_0=0$.

Suppose that a sequence $\{u_i\, :\, i\le 3j\}$ is already constructed for some $j\ge 0$. Then put
$$u_{3j+1}:=u_{3j}+B(3j-1)-C(3j-1),\, u_{3j+3}:=u_{3j+1}+A(3j+1)-C(3j+1),$$
\noindent while $u_{3j+2}$ define arbitrarily fulfilling the inequality $u_{3j+2}\ge u_{3j}+A(3j)-C(3j)$. Employing (\ref{1}) one can verify by induction on $j$ that the constructed sequence satisfies vector  $(A,\, B,\, C)$. To this end, it suffices to show that $u_{3j+2}\ge u_{3j+1}+A(3j+1)-B(3j+1),\, u_{3j+2}\ge u_{3j+3}+B(3j+2)-A(3j+2)$. 

%Hence, the entropy $H\ge 1/3$ since the coordinates $u_{3j+2},\, j\ge 0$ attain arbitrarily sufficiently large values. 
Therefore,
denoting by $U_N:=\{u_i\, :\, 0\le i<N\} \subset \RR^N$  the tropical prevariety of sequences satisfying (\ref{0}) for $0\le j<N-2$, we get $\dim (U_N)\ge \lfloor N/3 \rfloor$ since the coordinates $u_{3j+2},\, j\ge 0$ attain arbitrarily sufficiently large values. Hence, 
\begin{eqnarray}\label{19}
\lim \inf_{N\to \infty} \dim(U_N)/N\ge 1/3
\end{eqnarray}
 \vspace{2mm}

To prove the opposite inequality $H\le 1/3$ fix for the time being $N$, a sequence $(u_0,\dots,u_{N-1})\in U_N$, and consider a graph $G:=G_N$ with $N$ vertices $\{0,\dots,N-1\}$. For each $0\le j\le N-3$, if the minimum in (\ref{0}) is attained on some two indices among $j,\, j+1,\, j+2$ then between these indices draw an edge in $G$ (it is not excluded that there could be three edges between  $j,\, j+1,\, j+2$). Observe that any such graph $G$ determines a convex polyhedron in $U_N$ whose dimension does not exceed the number $c$ of connected components of $G$. Thus, our next goal is to bound $c$ from above.

We call an interval of $G$ a maximal (with respect to inclusion) subset (of at least two elements) of its vertices being an interval and belonging to the same connected component. Enumerate the connected components in an arbitrary way.  

\begin{lemma}\label{graph}
The following statements are valid under the assumption of either (\ref{1}) or (\ref{2}):

i) two intervals can't adjoin;

ii) the vertices between two neighbouring intervals belong to two alternating components $s,\, t$, respectively;

iii) if neighbouring intervals belong to components $p,\, q$, respectively, then $\{p,\, q\} \subset \{s,\, t\}$.
\end{lemma}   

{\bf Proof}. i) Suppose that $i,\, i+1$ belong to the same connected component, while $i+2,\, i+3$ belong to a different connected component. Then (considering triples of vertices $i,\, i+1,\, i+2$ and $i+1,\, i+2,\, i+3$) (\ref{0}) implies that 
$$(u_i+A(i)=)u_{i+1}+B(i)<u_{i+2}+C(i),$$
$$ u_{i+1}+A(i+1)>u_{i+2}+B(i+1)(=u_{i+3}+C(i+1))$$
\noindent which leads to a contradiction with (\ref{1}) as well as with (\ref{2}).

ii) If $i,\, i+1$ belong to the same connected component $p$, while $i+2$ belongs to a different component $s$ then $i+3$ belongs to $p$ (considering triple $i,\, i+1,\, i+2$). Similarly, $i+4$ belongs to $s$ (provided that $i+3$ does not belong to the next interval). Continuing this argument we establish ii) and in addition, iii) in case $p=q$.
To establish iii) in case $p\neq q$ note that $s=q$ considering triple $j-2,\, j-1,\, j$ where $j$ is the beginning of the next interval (belonging to $q$). $\Box$ \vspace{2mm}

One can deduce  statement 1) of the Theorem from Lemma~\ref{graph}. Indeed, for any new emerging component (by scanning $G$ from the left to the right) its first element is located betweem two neighbouring intervals because of i), in addition, between two neighbouring intervals at most one new component can emerge due to ii), iii).  This entails an upper bound $\lceil N/3 \rceil +1$ on the number of connected components of $G$ since the length of every interval is at least $2$. Hence the entropy $H\le 1/3$ and the limit $\lim_{N\to \infty} \dim(U_N)/N=1/3$ taking into account (\ref{19}).

\subsection{Case $A(i)+C(i-1)<B(i-1)+B(i),$
$$A(i+1)-A(i)+C(i)-C(i-1)=B(i+1)-B(i-1)$$}

Now we proceed to the proof of statement 2) of the Theorem. First similar to section~\ref{1.1} construct by induction a family of holonomic sequences $V_N:=\{v_i\, : \, 0\le i<N\}\subset \RR^N$ satisfying (\ref{0}) such that $\dim(V_N)\ge N/4 - const$.
Let $4j_0$ be greater than all the roots of polynomial $B(i-1)+B(i)-A(i)-C(i-1)$. Obviously, any sequence $v_0,\dots,v_i$ satisfying (\ref{0}) one can continue to $v_0,\dots,v_i,\, v_{i+1}$ also satisfying (\ref{0}). Construct an arbitrary sequence $v_0,\dots,v_{4j_0}$ satisfying (\ref{0}). Then family $U_{4j_0+1}$ consists of this single sequence. Consider the latter as a base of induction.

Suppose that a family $V_{4j+1}=\{v_i\, :\, 0\le i\le 4j\}$ is already constructed for some $j\ge j_0$. Put
$$v_{4j+1}:=v_{4j}+B(4j-1)-C(4j-1),$$
$$ v_{4j+2}:=v_{4j}+A(4j)-C(4j),$$
$$ v_{4j+4}:=v_{4j+2}+A(4j+2)-C(4j+2).$$
\noindent Then we take $v_{4j+3}$ in an arbitrary way satisfying inequality $v_{4j+3}\ge v_{4j+2}+B(4j+1)-C(4j+1)$. One can verify that the constructed family $V_{4j+5}:=\{v_i\, :\, 0\le i\le 4j+4\}$ satisfies (\ref{0}). To this end it suffices to show that 
$$u_{4j}+A(4j)<u_{4j+1}+B(4j),\, u_{4j+1}+A(4j+1)=u_{4j+2}+B(4j+1),$$
$$v_{4j+3}\ge v_{4j+2}+A(4j+2)-B(4j+2),\, v_{4j+3}\ge v_{4j+4}+B(4j+3)-A(4j+3)$$
\noindent employing (\ref{2}), (\ref{3}). Similar to  section~\ref{1.1} conclude that $\dim(V_{4j+1})=j-j_0$, thus 
\begin{eqnarray}\label{20}
\lim \inf_{N\to \infty} \dim(V_N)/N \ge 1/4
\end{eqnarray}
%the entropy $H\ge 1/4$. 
\vspace{2mm}

Now we prove the opposite inequality $H\le 1/4$ in statement 2) of the Theorem. It goes analogously to the proof of the upper bound on the entropy in section~\ref{1.1} considering graph $G$ and applying Lemma~\ref{graph} with a difference  that now the length of every interval is at least 3 (with a possible finite number of intervals of length 2 lying in $\{0,\dots,4j_0\}$).
Indeed, if $i,\, i+1$ was an interval then 
$$u_i+B(i-1)=u_{i+1}+A(i-1),\, u_i+A(i)=u_{i+1}+B(i),$$
\noindent and we get a contradiction with (\ref{2}). Thus, $H\le 1/4$ and (\ref{20}) implies the existence of the limit
$\lim_{N\to \infty} \dim(V_N)/N=H=1/4$. The
statement 2) of the Theorem is proved.

\section{Case of zero entropy} 

Now we proceed to a proof of statement 3) of the Theorem. 

\subsection{Case $A(i)+C(i-1)>B(i)+B(i-1)$}

First we assume that
\begin{eqnarray}\label{4}
A(i)+C(i-1)>B(i)+B(i-1)
\end{eqnarray}
(cf. (\ref{1}), (\ref{2})). Denote by $W_N:=\{w_i\, :\, 0\le i<N\}\subset \RR^N$ the set of all sequences satisfying (\ref{0}). We show that $\dim (W_N)$ is bounded from above by a constant independent from $N$.

\begin{lemma}\label{one}
Let $w:=\{w_i\, :\, 0\le i\le 4j_0\}\in W_{4j_0+1}$ (cf. section~\ref{1.1}). Then $w$ has at most one-dimensional continuation in any $W_N$ for $N>4j_0+1$.
\end{lemma}

{\bf Proof}. Let $\{w_i\, :\, 0\le i\le N\}\in W_{N+1}$ be a continuation of $w$. denote by $s\ge 4j_0$ the maximal  integer such that for every $4j_0\le t<s$ it holds $w_t+A(t)=w_{t+1}+B(t)$. Suppose w.l.o.g. that $s<N$. We claim that $\{w_i\, :\, 0\le i\le s+1\}\in W_{s+2}$ has a unique continuation in any $W_N$ for $N>s+1$.

We consider two cases. In the first one 
\begin{eqnarray}\label{5}
w_s+A(s)>w_{s+1}+B(s)
\end{eqnarray}
Then $w_{s+2}+C(s)=w_{s+1}+B(s)$ because of (\ref{0}). This together with (\ref{4}) implies that $w_{s+1}+A(s+1)>w_{s+2}+B(s+1)$. Continuing this arguing by induction on $i$ we obtain that $w_{i+2}+C(i)=w_{i+1}+B(i)$ for any $i\ge s$. This proves the claim in case (\ref{5}).

In the second case 
\begin{eqnarray}\label{6}
w_s+A(s)<w_{s+1}+B(s)
\end{eqnarray}
Then we have 
\begin{eqnarray}\label{7}
w_{s+2}+C(s)=w_s+A(s)
\end{eqnarray}
Therefore, $w_{s+1}+A(s+1)>w_{s+2}+B(s+2)$ since otherwise, summing up the opposite inequality with (\ref{6}), (\ref{7}),
we get a contradiction with (\ref{4}). Hence, we arrive to already considered in the first case inequality (\ref{5}) which proves the claim. \vspace{2mm}

Thus, all the continuations of $w$ in any $W_N$  for $N>4j_0$ have the following form. For some $4j_0\le s<N$ for each $4j_0\le t<s$ it holds $w_t+A(t)=w_{t+1}+B(t)$. After that, $w_{s+1}$ can take an arbitrary value not less than $w_{s-1}+A(s-1)-C(s-1)$ (when $s>4j_0$ or when $s=4j_0$ and $w_{4j_0-1}+A(4j_0-1)=w_{4j_0}+B(4j_0-1)$). Subsequently, a continuation of $\{w_i\, :\, 0\le i\le s+1\}$ in any $W_N$ for $N>s$ is unique. Lemma is proved. $\Box$ \vspace{2mm}

Lemma~\ref{one} implies that the entropy $H=0$ in case (\ref{4}).

\subsection{Case $A(i)+C(i-1)<B(i)+B(i-1),$
$$A(i+1)-A(i)+C(i)-C(i-1)>B(i+1)-B(i-1)$$}\label{2.2}

Now we assume (\ref{2}) and
\begin{eqnarray}\label{8}
A(i+1)-A(i)+C(i)-C(i-1)>B(i+1)-B(i-1)
\end{eqnarray}

\begin{lemma}\label{alternating}
Under assumptions (\ref{2}), (\ref{8}) let a sequence $\{w_i\, :\, 0\le i<N\}\in W_N$. If $w_{k+2}+C(k)=w_{k+1}+B(k)$ (cf. (\ref{0})) for a sufficiently large $k$ (namely, $k>4j_0$, cf. section~\ref{1.1}) then 
\begin{eqnarray}\label{9}
w_{k+2j+2}+C(k+2j)=w_{k+2j+1}+B(k+2j),
\end{eqnarray}
\begin{eqnarray}\label{10}
w_{k+2j+3}+C(k+2j+1)=w_{k+2j+1}+A(k+2j+1)
\end{eqnarray} 
for any $j\ge 0$.
\end{lemma}

{\bf Proof} goes by induction on $j$ . Suppose that (\ref{9}) holds for some $j\ge 0$ (the base is for $j=0$). Then $w_{k+2j+2}+B(k+2j+1)>w_{k+2j+1}+A(k+2j+1)$ because of (\ref{2}), therefore, (\ref{10}) is true (invoking (\ref{0})). Summing up
(\ref{9}), (\ref{10}) and employing (\ref{8}) for $i:=k+2j$, we obtain that
$$w_{k+2j+3}+B(k+2j+2)<w_{k+2j+2}+A(k+2j+2).$$
\noindent Hence, (\ref{9}) is valid for $j:=j+1$, this proves the inductive hypothesis and lemma. $\Box$ \vspace{2mm}

Consider graph $G:=G_N$ constructed in section~\ref{1.1} corresponding to the sequence $\{w_i\, :\, 0\le i<N\}$. Lemma~\ref{alternating} entails that $G$ contains at most one interval of the form $k,\, k+1, \dots$ where $k>4j_0$. Moreover, if there is such an interval then it is of the form  $k,\, k+1, \dots, N$, so ends at $N$. Therefore, $G_N$ can contain at most one new connected component in comparison with $G_{4j_0}$. Indeed, if $m>4j_0$ is the first element  of a new component then $m-2,\, m-1$ belong  to the same component (due to (\ref{0})), hence  $m-2,\, m-1$ is the end of an interval which contradicts to Lemma~\ref{alternating}. Note that $G_N$ can contain even less components than $G_{4j_0}$ does.

From this follows that the number of components of $G$ is bounded from above by a constant independent from $N$. Hence (cf. section~\ref{1.1}), $\dim(W_N)$ is bounded from above by a constant independent from $N$, thus the entropy $H=0$ under assumptions (\ref{2}), (\ref{8}).

\subsection{Case $A(i)+C(i-1)<B(i)+B(i-1),$
$$A(i+1)-A(i)+C(i)-C(i-1)<B(i+1)-B(i-1)$$}

\begin{lemma}\label{even}
Let (\ref{2}) and
\begin{eqnarray}\label{11}
A(i+1)-A(i)+C(i)-C(i-1)<B(i+1)-B(i-1)
\end{eqnarray}
be fulfilled. Assume that a sequence $\{w_i\, :\, 0\le i<N\} \in W_N$ and for some $k>4j_0$ it holds
\begin{eqnarray}\label{12}
w_{k+2}+C(k)=w_{k+1}+B(k).
\end{eqnarray}
Then for any $j\ge 1$ we have 
\begin{eqnarray}\label{13}
w_{k+j+2}+C(k+j)=w_{k+j}+A(k+j).
\end{eqnarray}
\end{lemma}

{\bf Proof} goes by induction on $j$. Suppose that (\ref{13}) is already established for all $1\le j\le s$ for some $s\ge 0$ (thus, one can treat the void case $s=0$ as a base of induction). 

%First let $s$ be odd. 
Summing up inequalities (\ref{11}) for $i:=k+2t,\, 0\le t\le r$ we obtain that
\begin{eqnarray}\label{14}
B(k)+\sum_{1\le m\le 2r} (-1)^m A(k+m)< B(k+2r)+\sum_{0\le m\le 2r-1} (-1)^m C(k+m)
\end{eqnarray}

First, let $s$ be odd.
Summing up equations (\ref{13}) multiplied by signs $(-1)^j$ for $1\le j\le s$, with (\ref{14}) for $r:=(s+1)/2$ and also with (\ref{12}) we get inequality 
\begin{eqnarray}\label{16}
w_{k+s+2}+B(k+s+1)>w_{k+s+1}+A(k+s+1).
\end{eqnarray}
\noindent Hence, (\ref{0}) entails that
$$w_{k+s+3}+C(k+s+1)=w_{k+s+1}+A(k+s+1)$$
\noindent which proves (\ref{13}) for $j:=s+1$.

Now let $s$ be even. Applying (\ref{14}) for $k:=k+1,\, r:=s/2$ we obtain inequality
$$B(k+1)+\sum_{1\le m\le s} (-1)^m A(k+m+1)< B(k+s+1)+\sum_{0\le m\le s-1} (-1)^m C(k+m+1).$$
\noindent Summing up the latter inequality with (\ref{2}) for $i:=k+1$, we establish inequality
\begin{eqnarray}\label{15}
-B(k)+\sum_{0\le m\le s} (-1)^m A(k+m+1)< B(k+s+1)+\sum_{0\le m\le s} (-1)^m C(k+m)
\end{eqnarray}
Similar to the case of odd $s$ summing up equations (\ref{13}) multiplied by signs $(-1)^j$ for $1\le j\le s$, with (\ref{15}) and with (\ref{12}), we again get inequality (\ref{16}). Then again (\ref{0}) entails (\ref{13}) for $j:=s+1$ which completes the proof of the inductive hypothesis and Lemma. $\Box$ \vspace{2mm}

As in section~\ref{2.2} Lemma~\ref{even} implies that graph $G:=G_N$ can contain at most one interval of the form $k,\, k+1,\dots$ where $k>4j_0$. Also as in section~\ref{2.2} we conclude that the entropy $H=0$. This completes the proof of the Theorem. $\Box$ \vspace{2mm}

{\bf Acknowledgements}. The author is grateful to the grant RSF 16-11-10075 and
 to MCCME for inspiring atmosphere.

\end{document}